\documentclass[english,aps,preprint]{revtex4}
\usepackage[T1]{fontenc}
\usepackage[utf8]{luainputenc}
\usepackage{bm}
\usepackage{amsmath}
\usepackage{graphicx}

\makeatletter

\usepackage[caption=false,font=normalsize,labelfont=sf,textfont=sf]{subfig}

\usepackage[T1]{fontenc}
\usepackage[utf8]{luainputenc}
\usepackage{bm}
\usepackage{amsmath}
\usepackage{graphicx}

\begin{document}
%

\title{A Convex Functional for Image Denoising based on Patches with Constrained Overlaps and its vectorial application to Low Dose Differential Phase Tomography }

%
%
%

\author{Alessandro Mirone}

\affiliation{European Synchrotron Radiation Facility, BP 220, F-38043 Grenoble
Cedex, France}

\author{Emmanuel Brun}

\affiliation{ European Synchrotron Radiation Facility and the 
 Ludwig Maximilians University, Physics Departement,  85748 Garching, Germany}

\author{Paola Coan}

\affiliation{Ludwig Maximilians University, Departement of Clinical Radiology,80336 Munich, Germany  
  and  LMU, Physics Departement,  85748 Garching, Germany}

\begin{abstract}
We solve the image denoising problem with a dictionary learning technique by writing a convex functional of a new form. This functional contains beside the usual sparsity inducing term and fidelity term, a new term which induces similarity between overlapping patches in the overlap regions. The functional depends on two free regularization parameters: a coefficient multiplying the sparsity-inducing $L_{1}$ norm of the patch basis functions coefficients, and a coefficient multiplying the $L_{2}$ norm of the differences between patches in the overlapping regions. The solution is found by applying the iterative proximal gradient descent method with FISTA acceleration. In the case of tomography reconstruction we calculate the gradient by applying projection of the solution and its error backprojection at each iterative step. We study the quality of the solution, as a function of the regularization parameters and noise, on synthetic datas for which the solution is a-priori known. We apply the method on experimental data in the case of Differential Phase Tomography. For this case we use an original approach which consists in using vectorial patches, each patch having two components: one per each gradient component. The resulting algorithm, implemented in the ESRF tomography reconstruction code PyHST, results to be robust, efficient, and well adapted to strongly reduce the required dose and the number of projections in medical tomography. 
\end{abstract}

\maketitle

\section{Introduction}

 Users of X-ray tomography aim to push the frontiers of their studies towards new domains which require finer time resolution, better signal to noise ratio, and less radiation damage. All these three requirements bring to a data-starving situation where for a given quality goal, the available data volume is never enough. A solution to this problem consists in filling in the gap left by the missing data with an a-priori knowledge of the solution. The signals occurring in Nature, when cleaned from noise, present most of the time an intrinsic sparsity when expressed in the proper basis. An image is intrinsically sparse when it can be approximated as a linear combination of a small number $n$ of basis functions with $n\ll N$, where $N$ is the image dimensionality. Piece-wise constant images are examples of sparse signal : they have non-zero signal only at the flat regions borders when they are expressed by their gradient. For piece-wise constant images one can apply very efficient methods based on minimization of a convex functional, said also convex objective function, which contains a total variation penalty term. For other classes
 of images such as medical images one has to choose different solutions which are adapted to the intrinsic sparsity of the case under study. There are mainly two ways : either the sparsity structure is a-priori known and an appropriate basis of functions can be built from the beginning, or it must be automatically learned from a learning set with the dictionary learning technique. The dictionary learning technique has recently been applied to tomography reconstruction using the Orthogonal Matching Pursuit (OMP) denoising procedure\cite{dual_dictionary}. This procedure consists in obtaining first an over-complete basis of functions and then in fitting every patch of the image using at most $N_{omp}$  components selected from this basis. The components are heuristically selected choosing, each time, the one having the maximum overlap with the remaining error. This optimization method cannot be implemented as a convex objective function optimization problem because the linear combination of two candidate solutions can have more than $N_{omp}$ components. In other words the optimization domain is not convex. In this paper we present a more advanced formalism based on a convex functional that we describe in section \ref{sec:metodo}. For the solution of our functional minimization problem we have applied the recently developed tools taken from the field of convex optimization\cite{combettes:hal-00643807}. The result is a robust and efficient algorithm that we discuss and illustrate with synthetic and medical data in section \ref{sec:application}.

\section{methods}

\label{sec:metodo}In this section we introduce first the decomposition of an image into non-overlapping patches and the related objective function for denoising. Then we introduce our original formalism which ensures, using overlapping patches, a smooth transition at the patches borders, and finally we apply this formalism to image denoising and tomography reconstruction.

We denote by ${\bf 1}_{p}$ the indicator function of patch p, which is equal to 1 over the patch support ( typically an m{*}m square) and is zero otherwhere. For non-overlapping patches,  covering the whole domain, we have :

\begin{equation}
\sum_{p}\bm{1}_{p}(i)=1\:\forall i.
\end{equation}

where $i$ denotes the pixel position and can be thought as a two-dimensional vector. We are looking for the ideal solution ${\bf x}$ that we express by the vector, ${\bf w}$, of its coefficients in the basis of patch functions: 

\begin{equation}
x_{i}=\sum_{p}\bm{1}_{p}(i)\sum_{k}w_{kp}\varphi_{k}(i-r_{p}).\label{eq:wx}
\end{equation}

where the set $\left\{ \varphi_{k}\right\} $ is an over-complete basis of functions over the patch support; $r_{p}$ is the closest to the origin corner of the patch $p$, and $w_{kp}$ is the component $k,p$ of vector ${\bf w}$ which multiplies the basis function $\varphi_{k}$ in the  patch $p$. The denoising problem, given an image ${\bf y}$, consists in finding the minimum of a functional $F({\bf w})=f({\bf w})+g({\bf w})$ which is sum of two terms. The term $f({\bf w})=\left\Vert {\bf y}-{\bf x}\right\Vert _{2}^{2}$ links the solution to the the data ${\bf y}$. The other term, $g({\bf w})$, contains the a-priori knowledge about the solution. This way of breaking the functional in two terms has his roots in the Bayesian theorem. From a probabilistic point of view the denoising problem consists in finding, given a noisy image $\underline{{\bf y}}$ of an object, the most probable object ${\bf x}$ that can generate that image. We represent the object ${\bf x}$ through the patches coefficients ${\bf w}$. The Bayes theorem, applied to denoising, states that the conditional probability of ${\bf w}$ being the exact object given a measurement ${\bf y}$, is the product of the probability of ${\bf y}$ being the measure given the exact solution ${\bf w}$, times the a-priori probability of ${\bf w}$.

Assuming gaussian noise, the conditional probability of ${\bf y}$ being the measure given the exact solution ${\bf w}$ is $exp\left(-\left\Vert {\bf y}-{\bf x}\right\Vert _{2}^{2}/\left(2\sigma^{2}\right)\right)$, where ${\bf x}$ is expressed through the patches coefficients ${\bf w}$ by equation \ref{eq:wx}. The a-priori probability of ${\bf w}$ is written as $exp\left(-g({\bf w})/\left(2\sigma^{2}\right)\right)$. The a-priori knowledge that the solution is sparse in our patches basis can be expressed by using the sparsity-inducing $L_{1}$ penalization \cite{bach_sparsity_inducing}: $g({\bf w})=\beta\left\Vert {\bf w}\right\Vert _{1}$. The most probable solution ${\bf w^{\star}}$ is obtained by finding the $F({\bf w})$ minimum:

\begin{eqnarray}
{\bf w^{\star}} =& argmin_{{\bf w}}(f({\bf w})+g({\bf w})); \nonumber \\  f({\bf w}) =& \left\Vert {\bf y}-{\bf x}\right\Vert _{2}^{2};  \quad g({\bf w})=\beta\left\Vert {\bf w}\right\Vert _{1}.
\end{eqnarray}

The solution can be obtained  by using the iterative shrinkage thresholding algorithm\cite{combettes:hal-00643807} (ISTA) iterations:

\begin{equation}
{\bf w}_{n+1}=T_{\beta\gamma}({\bf w}_{n}-\gamma\nabla f({\bf w}_{n}));\quad{\bf w^{\star}}={\bf w}_{\infty}.
\end{equation}

where $T_{\alpha}$is the shrinkage operator defined as

\begin{equation}
T_{\alpha}({\bf w})=\frac{{\bf w}}{\left\Vert {\bf w}\right\Vert _{2}}max(\left\Vert {\bf w}\right\Vert _{2}-\alpha,0).
\end{equation}

and $\gamma$ is a positive number lesser than the inverse of the Lipschitz condition number $L$:

\begin{equation}
\gamma\in]0,1/L].
\end{equation}

The Lipschitz number L is such that :

\begin{equation}
\left\Vert \nabla f({\bf w}_{2})-\nabla f({\bf w}_{1})\right\Vert _{2}\leq L\left\Vert {\bf w}_{2}-{\bf w}_{1}\right\Vert _{2};\quad\forall{\bf w}_{1},{\bf w}_{2}.
\end{equation}

The ISTA algorithm can be accelerated by the Fast Iterative Shrinkage Thresholding \cite{Beck:2009:FIS:1658360.1658364}(FISTA) method.

In its non-overlapping version, the image denoising with patches is able to detect features that are within the field of the patch: if a line crosses the central region of a patch, it will be detected if the basis of function has been trained to detect such lines. But in the situation where a line intersects only one point in a corner of the patch square, the signal of this point is indistinguishable from that of a noisy point, no matters the dictionary training.

For this reason the patches denoising technique is often used with overlapping patches using post-process averaging\cite{elad_aharon}. In this case the minimization problem is solved for each patch separately first, and then the averaging is performed in the overlapped regions. 

In this study we do not follow this procedure but we add an overlap term into the objective function. We choose a system of patches which covers the whole domain, and we allow for overlapping. In this case the sum of all indicator functions is greater or equal to one :

\begin{equation}
\sum_{p}\bm{1}_{p}(i)\geq1;\:\forall i.
\end{equation}


We define the core indicator functions $\bm{1_{p}^{c}}$, which indicate the core of the patches, and make a non-overlapping covering:

\begin{equation}
\bm{1_{p}^{c}}(i)\leq\bm{1_{p}}(i);\quad\sum\bm{1_{p}^{c}}(i)=1;\:\forall i.
\end{equation}

For a given point $i$, $\bm{1_{p}^{c}}(i)$ indicates which patch
$p$ has its center $C_{p}$ closest to point $i$:

\begin{equation}
\sum_{p}\bm{1_{p}^{c}}(i)\left\Vert i-C_{p}\right\Vert _{1}\leq\left\Vert i-C_{p^{\prime}}\right\Vert _{1};\quad\forall p^{\prime},i.
\end{equation}

The solution ${\bf x}$ is composed using the central part of the
patches as indicated by the functions $\bm{1_{p}^{c}}$:

\begin{equation}
x_{i}=\sum_{p}\bm{1}_{p}^{c}(i)\sum_{k}w_{kp}\varphi_{k}(i-r_{p}).
\end{equation}

Now we introduce the $P$ operator which is the projection operator,
for tomography reconstruction, and is the identity for image denoising.
The functional $F({\bf w})$ whose minimum gives the optimal solution
is written, for both applications, as:

\begin{eqnarray}
F({\bf w})& =& f({\bf w})+g({\bf w});\quad g({\bf w})=\beta\left\Vert {\bf w}\right\Vert _{1};\nonumber\\
f({\bf w})& =& \left\Vert {\bf y}-{\bf P}({\bf x})\right\Vert _{2}^{2}+\nonumber\\
& & \rho\sum_{pi}\bm{1}_{p}(i)\left(x_{i}-\sum_{k}w_{kp}\varphi_{k}(i-r_{p})\right)^{2}. ~~
\end{eqnarray}

where the $\rho$ factor weights a similarity-inducing term which
pushes all the overlapping patches, which touch a point $i$, toward
the value $x_{i}$ of the global solution $x$ in that point. For future
reference we call $\left\Vert {\bf y}-{\bf P}({\bf x})\right\Vert _{2}^{2}$
the fidelity term.

The solution is found with the FISTA method, using the gradient of
$f({\bf w})$ which is easily written in compact form:

\begin{multline}
  \frac{\partial f({\bf w})}{\partial w_{kp}}=\sum_{i}2\varphi_{k}(i-r_{p})\bm{1}_{p}^{c}(i)\left\{ \left({\bf P}^{T}\left({\bf P}({\bf x})-{\bf y}\right)\right)_{i} 
 \begin{matrix} ~\\ ~ \\~ \end{matrix} \right. +\\
 \left.  \rho\sum_{p^{\prime}}\bm{1}_{p^{\prime}}(i)\left(x_{i}-\sum_{k^{\prime}}w_{k^{\prime}p^{\prime}}\varphi_{k^{\prime}}(i-r_{p^{\prime}})\right)\right\}+\\
  \sum_{i}2\varphi_{k}(i-r_{p})\rho{\bf 1}_{p}(i)\left(\sum_{k^{\prime}}w_{k^{\prime}p}\varphi_{k^{\prime}}(i-r_{p})-x_{i}\right)
\end{multline}

where $P^{T}$is the adjoint operator of P, and is called back-projection
operator in the case of tomography, and is still the identity for
image denoising.

\section{APPLICATION}

\label{sec:application} In this section we compare the fit with patches of a given image with or without overlapping. We analyze the quality of the recovered image as a function of the two regularization parameters $\beta$ and $\rho$. The same study is then performed on the same image with additive noise. Then we apply the method to tomography reconstruction of virtual phantoms using our convex functional with patch overlapping. We analyze the quality of the reconstruction as a function of the two regularization parameters projections. Finally, we apply the method onto experimental cases for the reduction of number of projection in order to reduce the deposited dose during a tomography.

\begin{figure}[htp]
 \subfloat{\includegraphics[width=0.45\columnwidth,height=0.45\columnwidth]{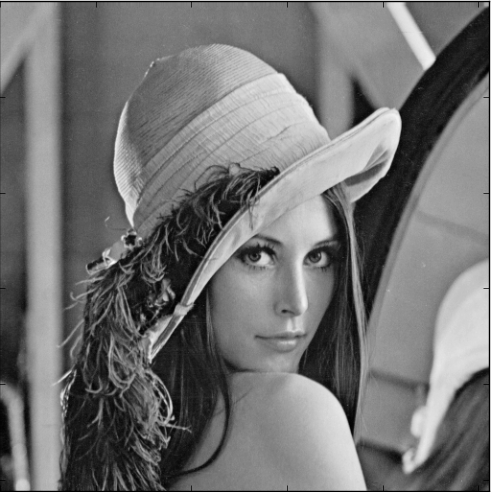}} 
  \hfill
 \subfloat{\includegraphics[width=0.45\columnwidth,height=0.45\columnwidth]{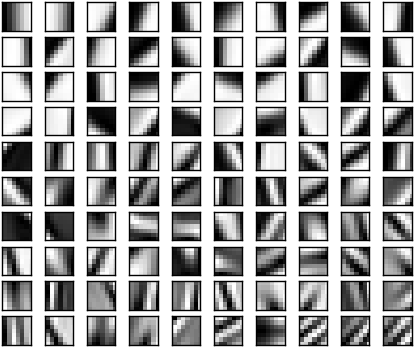}} 
 \caption{A training image and its K-SVD basis}
 \label{fig:A-training-image}
\end{figure}

\subsection{Fitting problem}

Figure \ref{fig:A-training-image} shows a training image (Lena) on the left. From this image we have extracted by dictionary learning the basis of patches that we show on the right. To obtain the basis we have implemented the K-SVD algorithm\cite{DBLP:conf/icpr/MazharG08} that we have applied using four atoms in the OMP procedure. Our basis consists in an over-complete set of 100 functions having a $7\times 7$ pixels support.

%

\begin{figure*}[htp]
\subfloat{\includegraphics[width=.22\columnwidth]{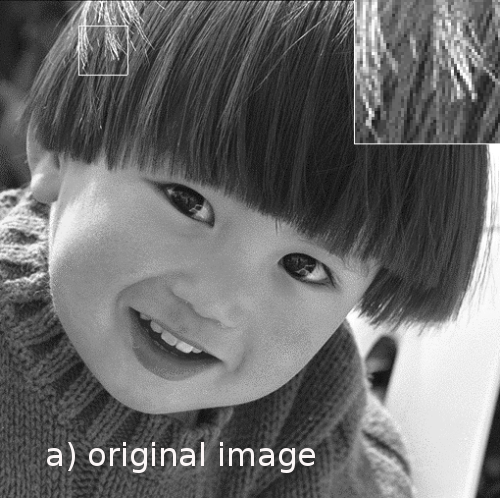}} 
\hfil
\subfloat{\includegraphics[width=.22\columnwidth]{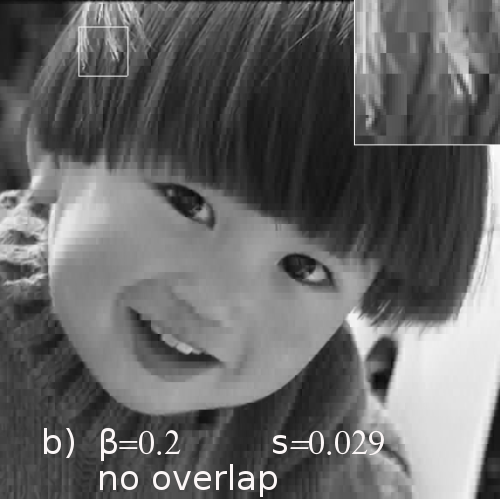}} 
\hfil
\subfloat{\includegraphics[width=.22\columnwidth]{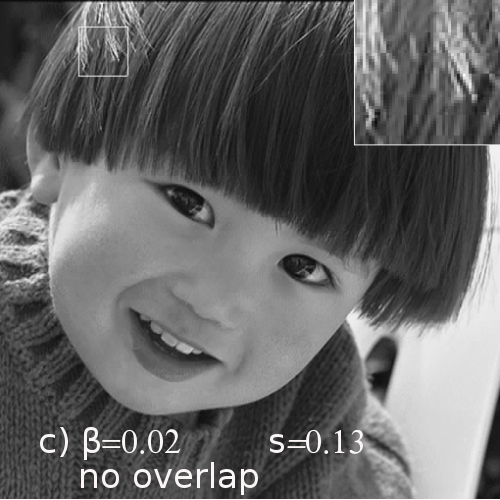}} 
\hfil
\subfloat{\includegraphics[width=.22\columnwidth]{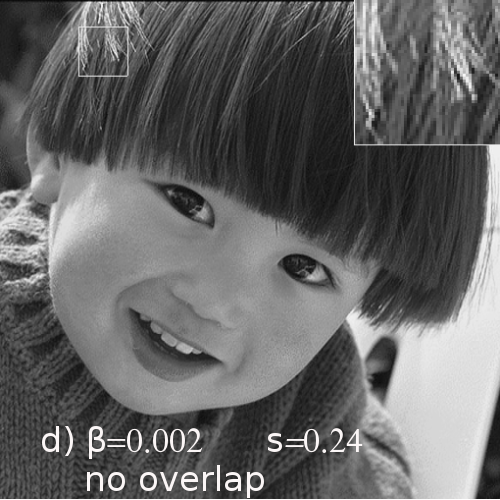}} 

\subfloat{\includegraphics[width=.22\columnwidth]{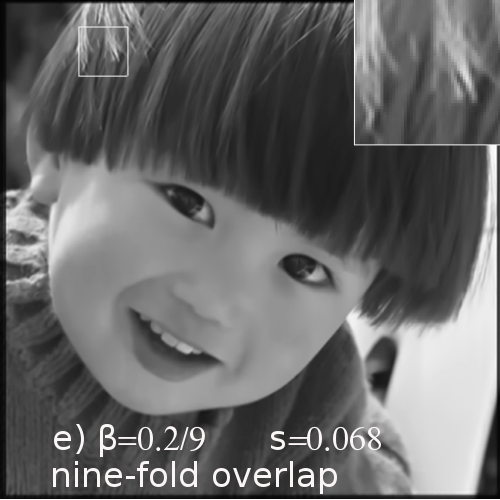}} 
\hfil
\subfloat{\includegraphics[width=.22\columnwidth]{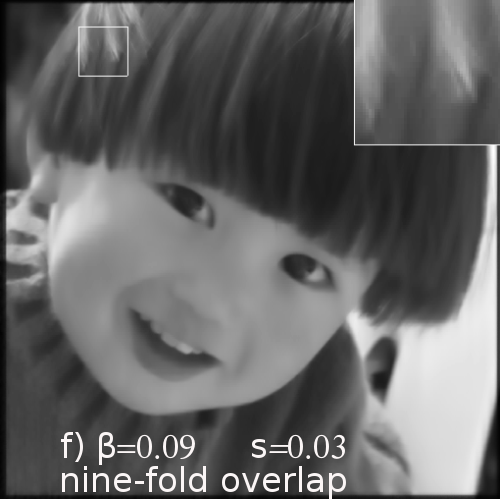}} 
\hfil
\subfloat{\includegraphics[width=.22\columnwidth]{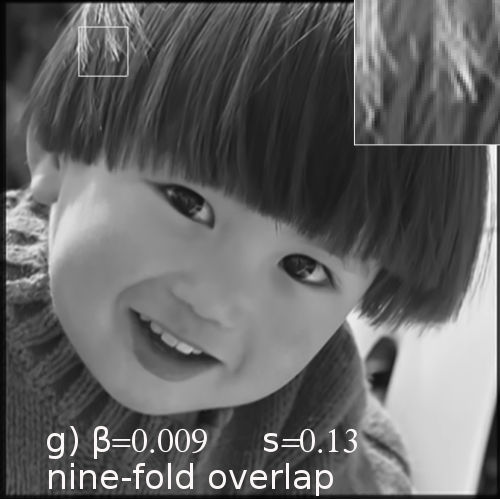}} 
\hfil
\subfloat{\includegraphics[width=.22\columnwidth]{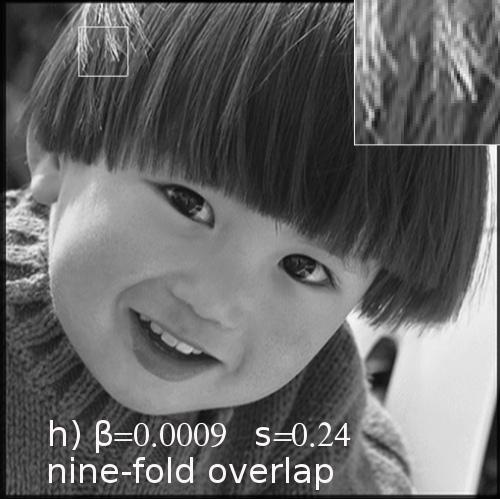}} 
\caption{ Fitting with no-overlapping (first row) and overlapping (second row) patches a noiseless image. The result for $\beta=0.2/9$ (e) reproduces the same fidelity/$L_{1}$ratio obtained with$\beta=0.2$ (b) in the non-overlapping case. For the same level of fitted details the overlapping method requires a higher number of components ( bigger $s$) in order to satisfy the overlapping constraints.} 
\label{fig:noov_noiseless}
\end{figure*}

We apply this set of patches that we have learned from Lena on a completely different image (Boy \cite{boy1}). In figure \ref{fig:noov_noiseless} we show the effect of fitting with non-overlapping patches a noiseless image. The original image is subfigure \ref{fig:noov_noiseless}a). It has been normalized to have its maximum equal to 1 and its minimum equal to 0. From left to right, of the first row, in the three columns on the right, we have the results for beta equal to 0.2, 0.02, 0.002. The sparsity $s$ is the ratio of non-zero components over the total number of components (100 components per every patch here). The value of sparsity for each image is reported on its legend. When we take the sparsity average over all the image, it takes the values 0.029,0.13,0.24 when beta goes from 0.2 to 0.002 in the first row of our example. 

From figure \ref{fig:noov_noiseless} it is obvious how the fit with non overlapping patches gives discontinuity at the patches borders. The upper-right corners of the figures are zooms in the hairs zone. We can see that, in the case of non-overlapping patches, features in form of lines are more sensitive to the regularization term when they are close to the patch border than when they are close to the center. When the regularization parameter is weakened by lowering beta, then more components are available for the fit and the dimensionality gets high enough to fit all the features. 

The second row of figure \ref{fig:noov_noiseless} shows the same fitting problem solved with our overlapping patches method. The set of overlapping patches has been obtained by replicating the original set of non-overlapping patches by translation vectors $(3l,3n)$; where $l$ and $n$ are integers in the interval $[0,2]$. This choice of translation vectors results in a nine-fold overlap : each point is covered by nine patches. The $\beta$ values take, in the last three columns of the second row, the values 0.09,0.009, 0.0009 from left to right. The overlap constraining factor is $\rho=1$. With these values these images have nearly the same sparsity than the non-overlapping results in the same columns of the first row. 

One can see in figure \ref{fig:noov_noiseless} that using overlapping patches the discontinuities have disappeared and the lines are homogeneously represented along their length even for the strongly regularized case. The overlap constraint reduces the tessellation effects but we can see that for the same sparsity ratio, another kind of distortion, in the form of smoothing, appears. This is due to the fact that the fidelity term, which links the data to the solution, is proportional to the area of the image, while the $L_{1}$ term acts on all the $w_{pk}$ coefficients whose number increases not only with the area of the image but also with the number of replica. Taking into account that we have 9 replica per patch, we show in the first image of second column the result obtained with $\beta=0.2/9$. This rescale back the $L_{1}$ term to a situation which is comparable with the $\beta=0.2$ case of image b) of the first row. The level of preserved details is comparable in these two cases. The sparsity, however is bigger in the overlapping case: $s=0.068$ which means nearly seven functions per patch. In the non-overlapping case, instead, it was about 3 components per patch. This means that the similarity inducing term is constraining a part of the available degrees to realize smooth transitions between neighboring patches.

\subsection{Denoising problem}

For the denoising problem two things must be kept in mind to select the regularization weights: for an image with a stronger noise a stronger value of $\beta$ will be necessary, but it is also true that for a stronger value of the regularizing parameter more features of the noiseless image will be filtered out. The best value of the regularization parameters is therefore a compromise between the necessity of filtering out the noise and that of preserving image features. 

\begin{figure}[htp]
 \subfloat{\includegraphics[width=0.22\columnwidth]{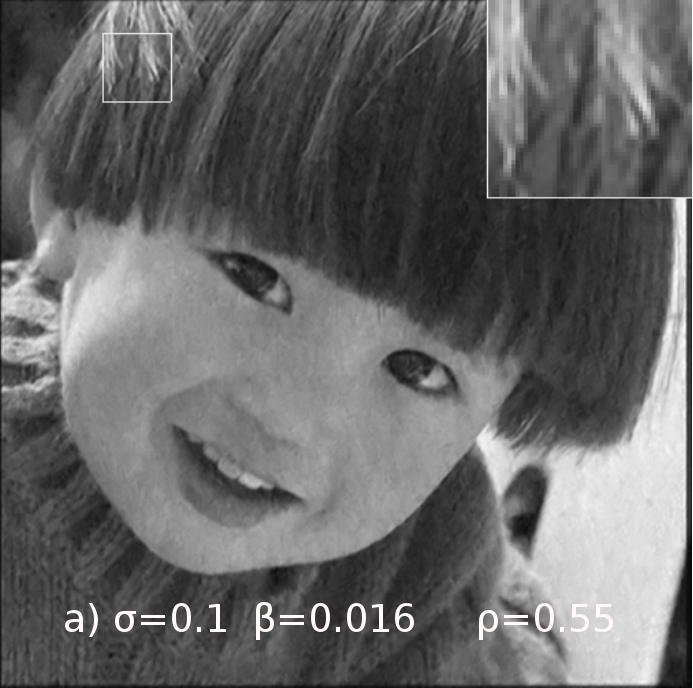}} 
 \hfill
 \subfloat{\includegraphics[width=0.22\columnwidth]{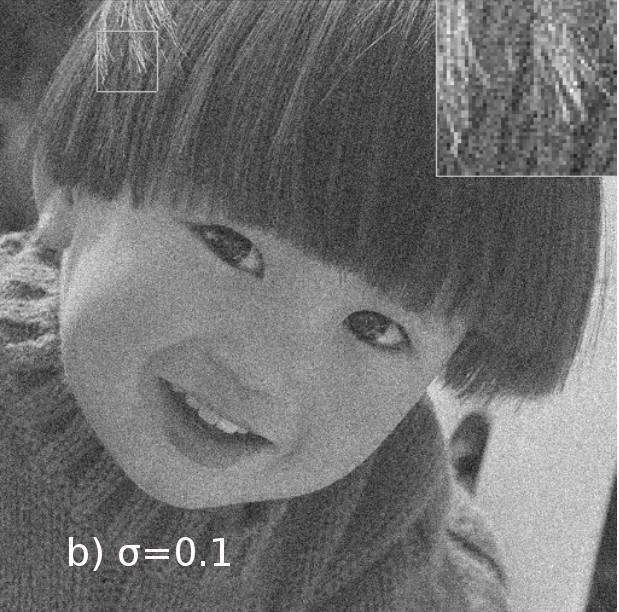}} 
 \hfill
 \subfloat{\includegraphics[width=0.22\columnwidth]{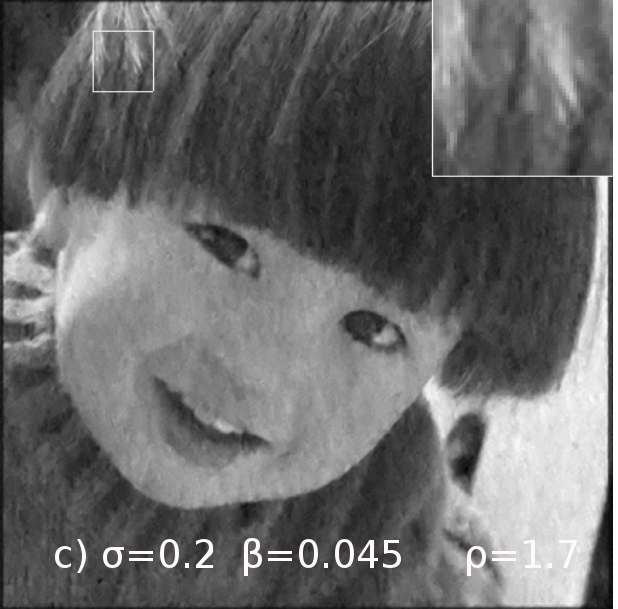}} 
 \hfill
 \subfloat{\includegraphics[width=0.22\columnwidth]{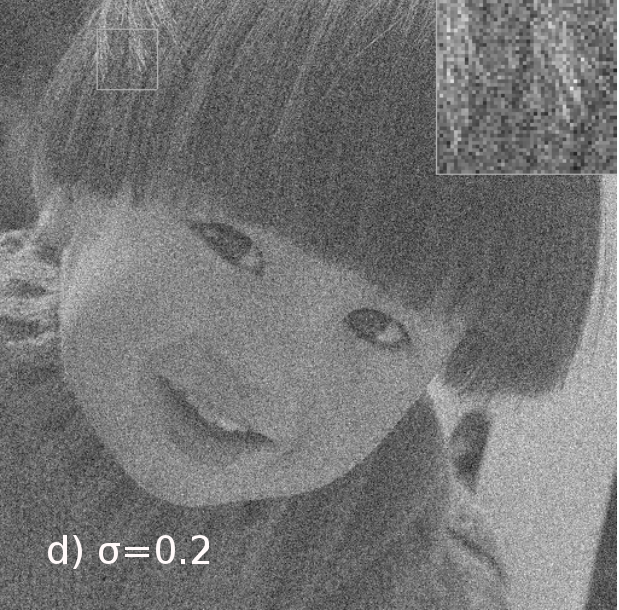}} 
 
 \caption{the effect of overlapping denoising (a,c) on an image with $\sigma=0.1,0.2$ (b,d) added white gaussian noise. The denoising has been performed with a replication step of 3 and the optimal $\beta,\rho$ 
.}
\label{fig:ov_noise}
\end{figure}

%
%
%
%


We define the quality improvement factor $Q$ based on the the Structural SIMilarity (SSIM \cite{ssim}) index S which is 1 when the images
are identical. Our definition of Q is : 
\begin{equation}
Q=\left(1-S({\bf y}_{noised},{\bf y})\right)/\left(1-S({\bf y}_{denoised},{\bf y})\right)\label{eq:qssim}
\end{equation}

\begin{figure} [htp]
 \subfloat{\includegraphics[width=0.99\columnwidth]{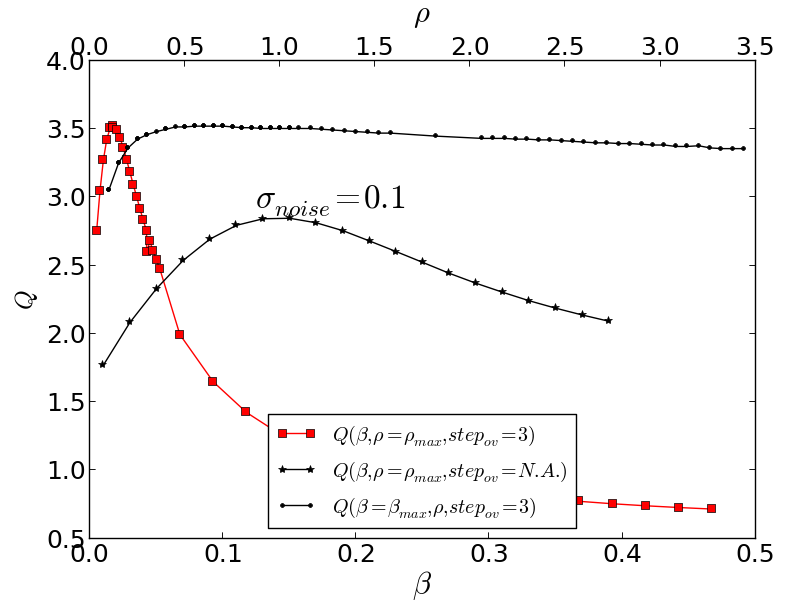}} 
 
 
 \caption{The quality improvement factor $Q$ for different values of $\beta,\rho$ around the optimal values for $\sigma_{\{noise\}}=0.1$. 
 For the  overlapping case a replication step of 3 has been used and we have varied separately $\beta$ (squares) and $\rho$(circles) around the optimal values. For the non-overlapping case only $\beta$ is varied.}
\label{fig:tendenze}
\end{figure}

We show in figure \ref{fig:tendenze} the quality improvement factor $Q$ for different values of $\beta,\rho$ for the overlapping case and the non-overlapping one. The tests have been done adding a white gaussian noise ($\sigma_{\{noise\}}=0.1$ and $0.2$) to the original image of figure \ref{fig:noov_noiseless}. For the overlapping case a replication step of 3 has been used and we have varied separately $\beta$ (squares) and $\rho$(circles) around the optimal values. In the non-overlapping case (stars) only $\beta$ is varied.

We observe that the quality improvement factor, for the overlapping case, peaks at a $\beta$ value which is about $1/9$ of the $\beta$ value of the no-overlapping peak. As explained before this lower $\beta$ is compensated by the nine-fold increase in the number of components which are weighted by $L_{1}$norm. We note that the optimal quality is better with the overlapping method. The trend of $\beta$ as a function of the noise level is as expected : a stronger noise needs to be regularized with a stronger $\beta$.

In figure \ref{fig:ov_noise} we show the denoising results (a and c) for two noised images ${\bf y}_{noised}$ (b and d) obtained adding a $\sigma_{\{noise\}}=0.1,0.2$ strong white gaussian to the original image : ${\bf y}_{noised}={\bf y}_{original}+noise(\sigma_{\{noise\}})$. The denoising has been performed with a replication step of 3 pixels and with the $\beta$ and $\rho$ giving the optimal SSIM. 

\subsection{Dose reduction in computed tomography}

\begin{figure}[htp]
 \subfloat{\includegraphics[width=0.49\columnwidth]{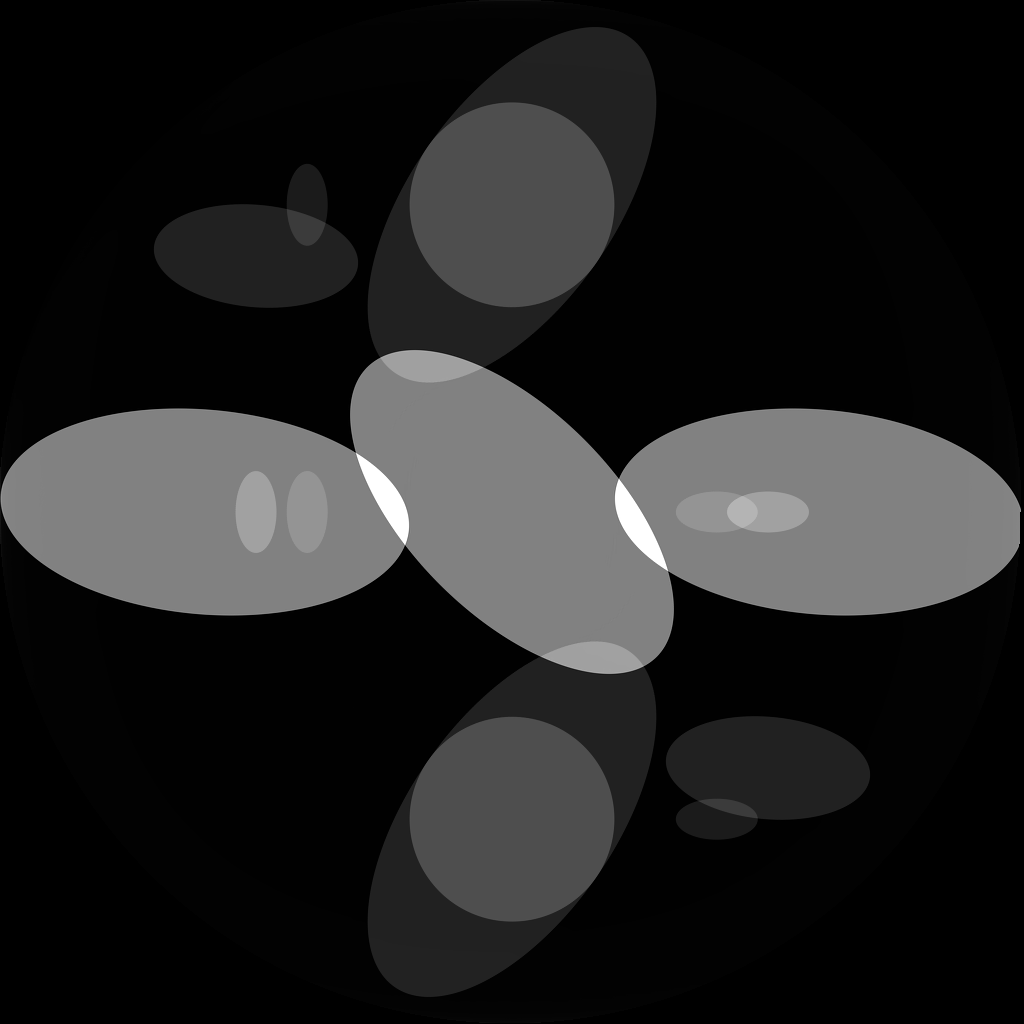}} 
  \hfill
 \subfloat{\includegraphics[width=0.49\columnwidth]{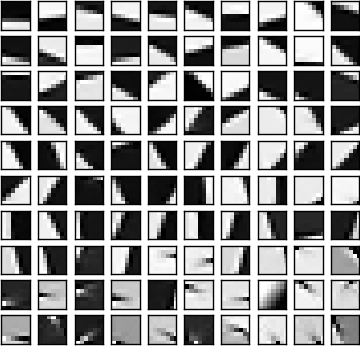}} 
 \caption{A phantom and its K-SVD basis        \label{fig:A-training-image-phantom}}
 
\end{figure}

In this section we apply the method using overlapping patches for tomography reconstruction on synthetic data. The figure \ref{fig:A-training-image-phantom} presents a $1024\times 1024$ pixels phantom and a bases of patches that has been learned from the phantom . 

\begin{figure}[htp]
\begin{center}
 \subfloat{\includegraphics[width=0.7\columnwidth]{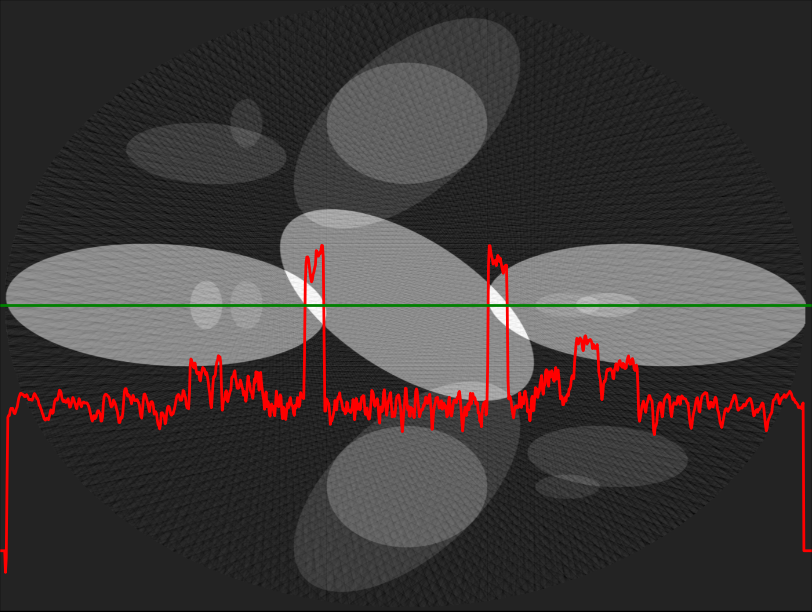}} 
  
 \subfloat{\includegraphics[width=0.7\columnwidth]{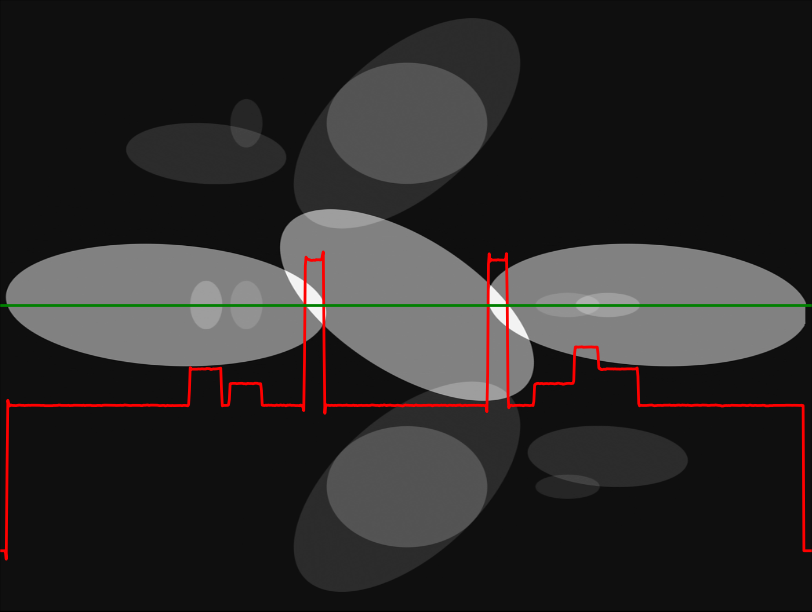}} 
 \caption{reconstruction from 150 projections of a $1024\times 1024$ phantom. Top: FBP; Down: our method with $\beta=0.001$ and $\rho=50$. \label{fig:pantom00fbppp}}
 \end{center}
\end{figure}

We show in figure \ref{fig:pantom00fbppp} the reconstruction obtained from 150 projections between 0 and 180 degree. Note that this is well below the number of projections, of the order of thousand, that would be required in this case by the Shannon-Nyquist sampling theorem. The lack of information due to the small number of projections results in a noised reconstruction when the filtered-back-projection reconstruction is applied (left) while the image recovered with our method (right), with $(3l,3n)$ translations for replication, maintains a high quality.

%
%
%
%
%
%

For this noiseless case the choice for $\beta$ and $\rho$ was easy and figure \ref{fig:pantom00fbppp} was produced with our first guess which gave already an excellent reconstruction, without need of further optimization. The first guess was $\beta=0.001$ and $\rho=50$. The quality improvement $Q$ between the FBP result and our method is $20$. The improvement factor is in this case the FBP error (the SSIM index distance from 1 as in \ref{eq:qssim}) divided by the error of our method.

%
%
%
%
%
%
%

\begin{figure}
\begin{center}
\includegraphics[width=0.95\columnwidth]{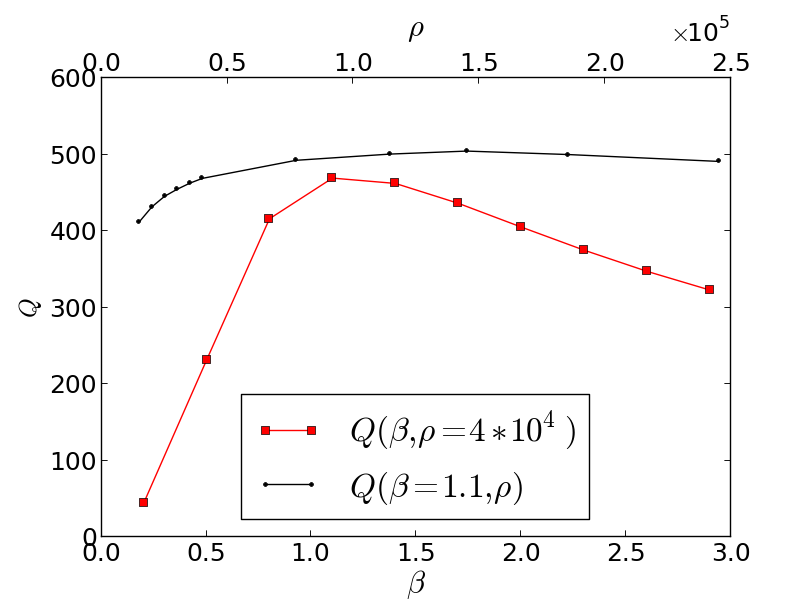}%
\caption{the $Q$ dependency versus $\beta$(squares) and $\rho$(circles) not far from the optimal choice of $\beta$ and $\rho$ for the reconstruction of the phantom from 150 noisy projection.}

\label{fig:opti-phantom} 
\end{center}
\end{figure}

The noisy case is instead more complex. To study it we add to the phantom sinogram a gaussian white noise with a $\sigma$ equal to $5\%$ of the maximum sinogram value. In figure \ref{fig:opti-phantom} we show the $Q$ dependency versus $\beta$(squares) and $\rho$(circles) not far from the optimal choice. We have not fully optimized $Q$ because we have found that in this extreme case where $\rho$ tends to be very big, we get a substantial decrease in convergence speed after the shoulder in the $\rho$(circles) curve. Before the shoulder we have convergence in some hundreds of iterations while after the shoulder it takes some thousands to converge.

We observe that at these high values of $\rho$ the effective behavior of our method is similar to those of a total variation penalization : the image is flattened over large regions. We think that the increase in convergence time is due to the fact that the algorithm takes a longer time to propagate when the flattened regions are larger.

\begin{figure}[htp]
\begin{center}
 \subfloat{\includegraphics[width=0.7\columnwidth]{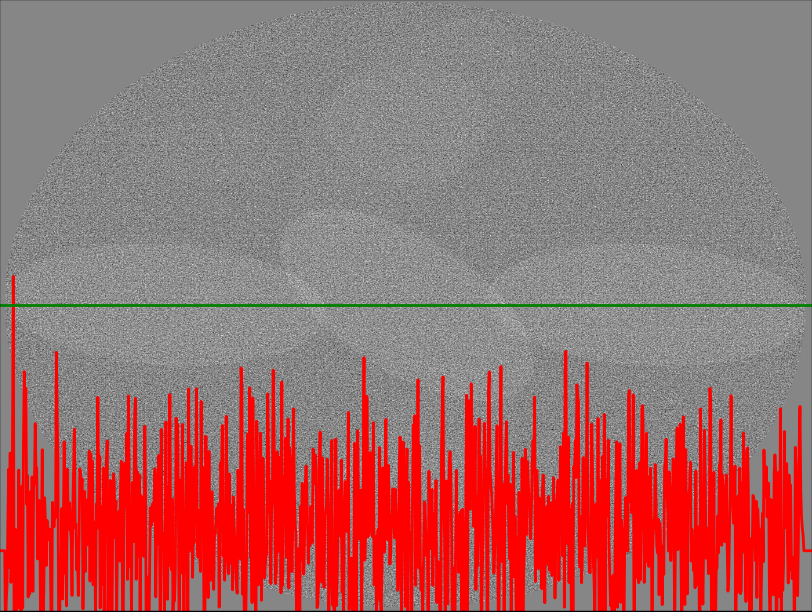}} 
  
 \subfloat{\includegraphics[width=0.7\columnwidth]{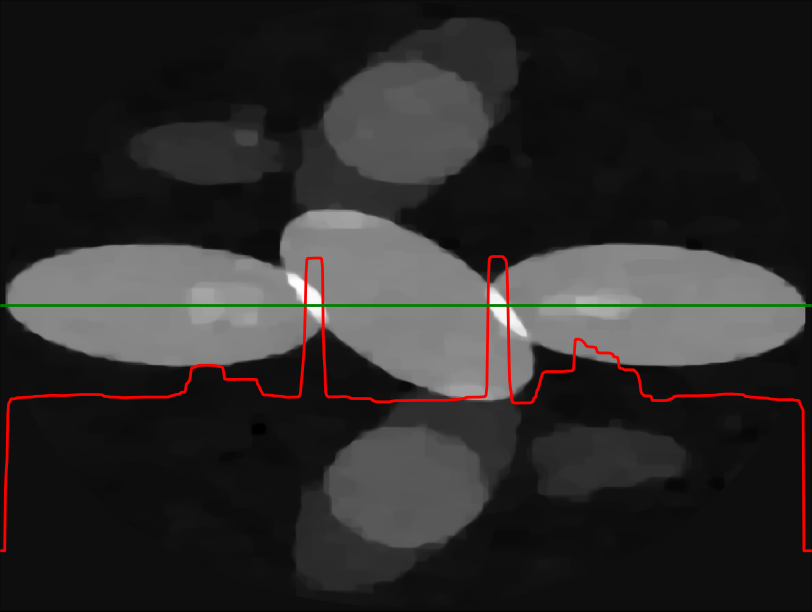}} 
 \caption{Reconstruction from 150 angles of the phantom sinogram after adding a gaussian white noise with a $\sigma$ equal to $5\%$ of the maximum sinogram value. Top : FBP. Bottom: our method with $\beta=1.1$ , $\rho=2.5*10^{5}$. \label{fig:pantom01fbppp}}
 \end{center}
 
\end{figure}

Figure \ref{fig:pantom01fbppp}, bottom, shows the result obtained with our method using the $\beta=1.1$ and $\rho=2.5*10^{5}$ values obtained from the $Q$ optimization illustrated figure \ref{fig:opti-phantom}. The FBP reconstruction is reported at the top of the figure. On these images it is clear that that our method outperforms standard FBP. The ellipsoids become visible and the plotted profile follows the expected behavior, whilst the FBP reconstruction shows only noise.
%
%
%
%
%

%
%
%

%
%
%
%
%
%
%

\subsection{Differential Phase tomography}

Finally we show a promising application of our method to medical tomography for a sample imaged using X-ray phase contrast imaging(PCI). This technique has shown an enhancement of soft tissue visualization in comparison to conventional techniques \cite{bravin2013}. PCI employs the dual property of X-rays of being simultaneously absorbed and refracted while passing through a tissue. Among all the phase contrast techniques we chose to test our method on the analyzer based imaging because  of the high sensitivity and unique results provided by this modality
for investigating large and highly absorbing biological tissues (i.e. full human breasts) \cite{Zhao2012}. Because breast are highly radio sensitive organs, X-ray CT of these organs are not clinically applied, even if a 3D would be benefit for radiologists. A reduction of the deposited radiation dose in CT combined with the unprecedented contrast improvement offered by PCI is thus of high interest for breast cancer detection.

In the analyzer based PCI technique, the projection data contain a signal which is proportional to the gradient of the X-ray phase in one direction (i.e. the direction perpendicular to the plane formed by the incoming and diffracted X-rays on a perfect Bragg crystal which is used for analyzing the radiation passing through the sample). When the object is rotated around an axis (Z-axis, for instance), this signal contains contributions from the X and Y gradient components, where the X and Y axis corotate with the sample. The two components are de-phased by a rotation angle of 90 degrees and can be reconstructed separately by multiplying before-hand the sinogram with the cosine and sine of the rotation angle. We apply our formalism considering that the reconstructed and learning images are vectorial objects : the value associated with a pixel is not a scalar but a two-component vector. More details on the principles and technical aspects of PCI are available in  \cite{bravin2013}. 

\begin{figure}[htp]
\begin{center}
\subfloat{\includegraphics[width=0.4\columnwidth]{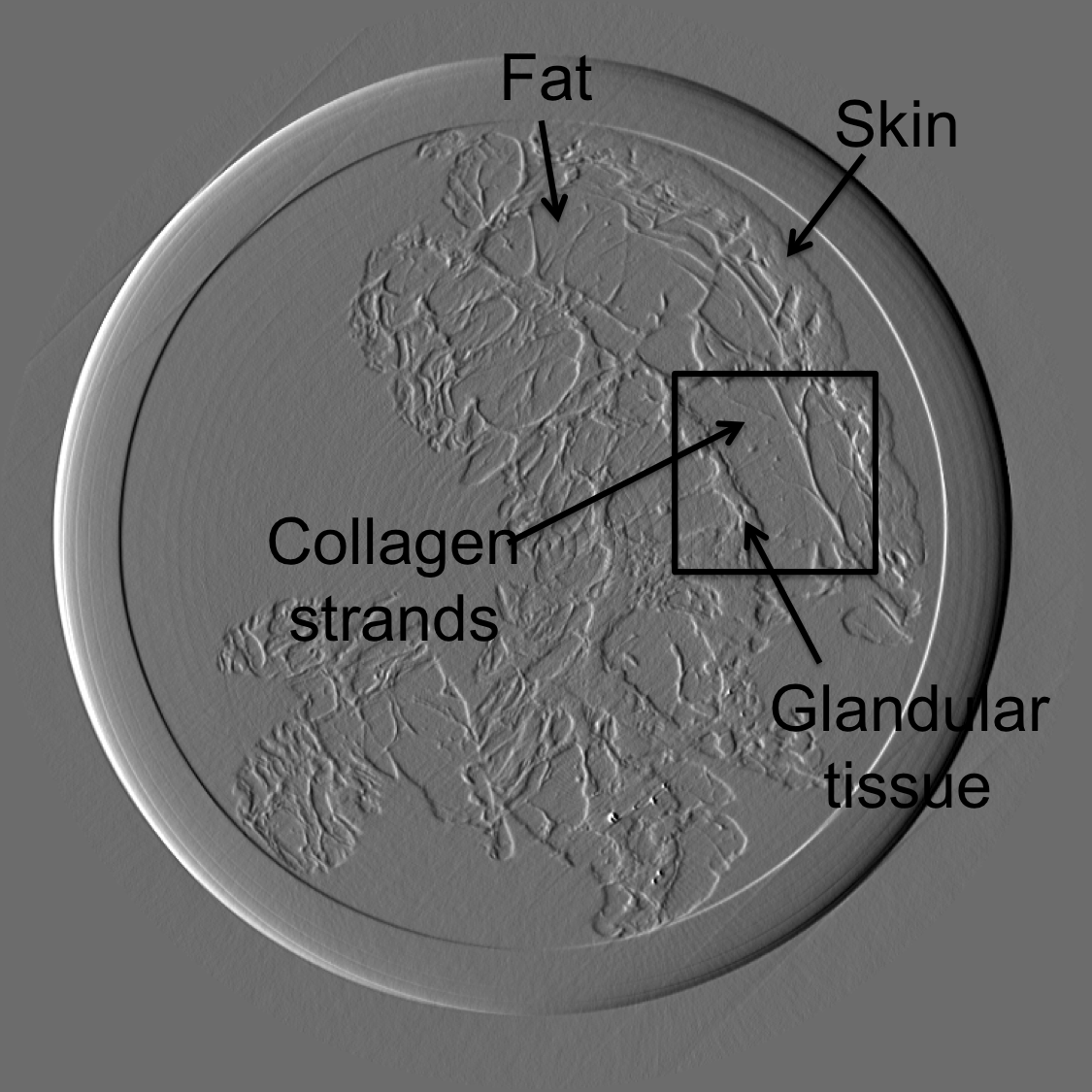}}

 \subfloat{\includegraphics[width=0.245\columnwidth]{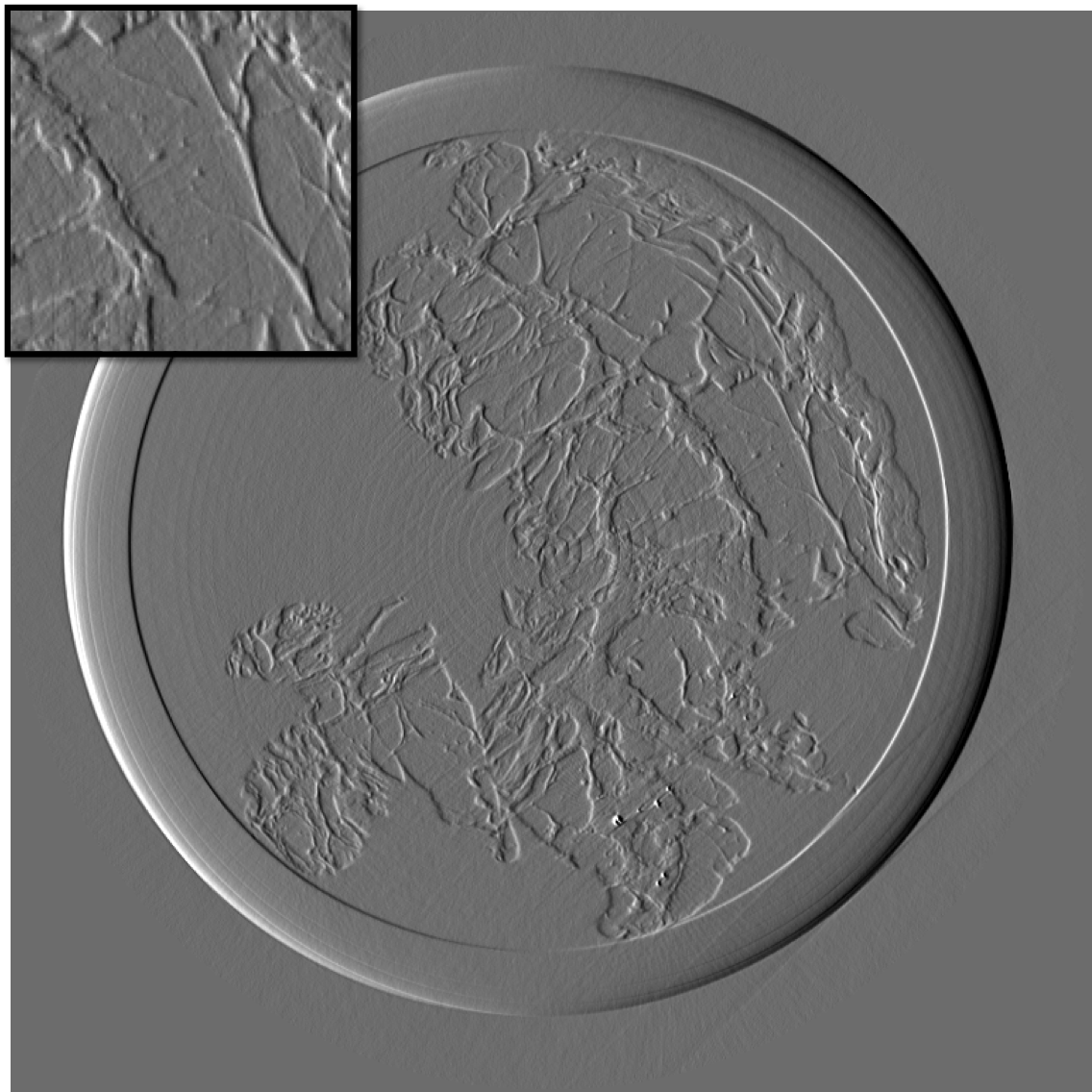}} 
  \hfill{}
   \subfloat{\includegraphics[width=0.245\columnwidth]{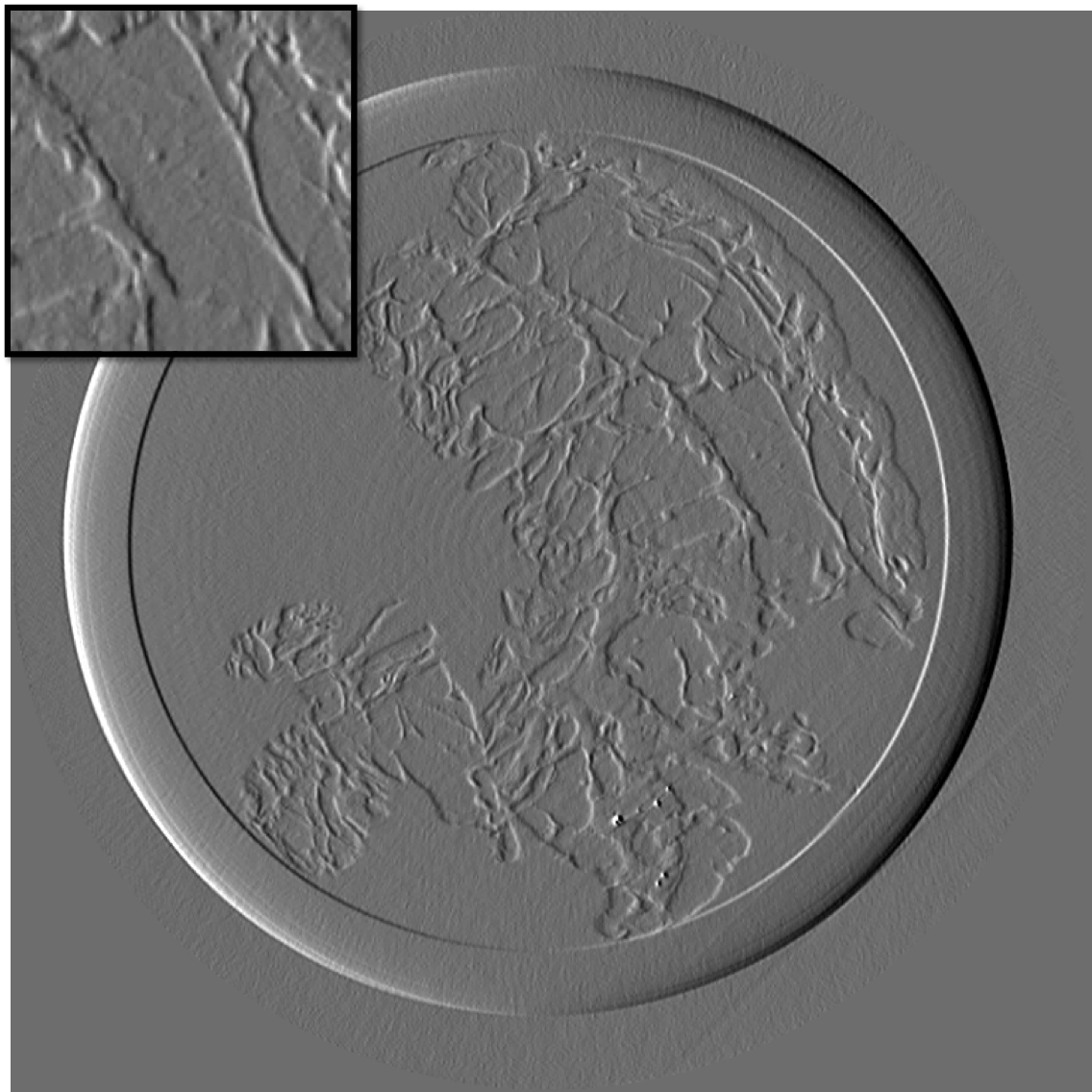}} 
 
 \subfloat{\includegraphics[width=0.245\columnwidth]{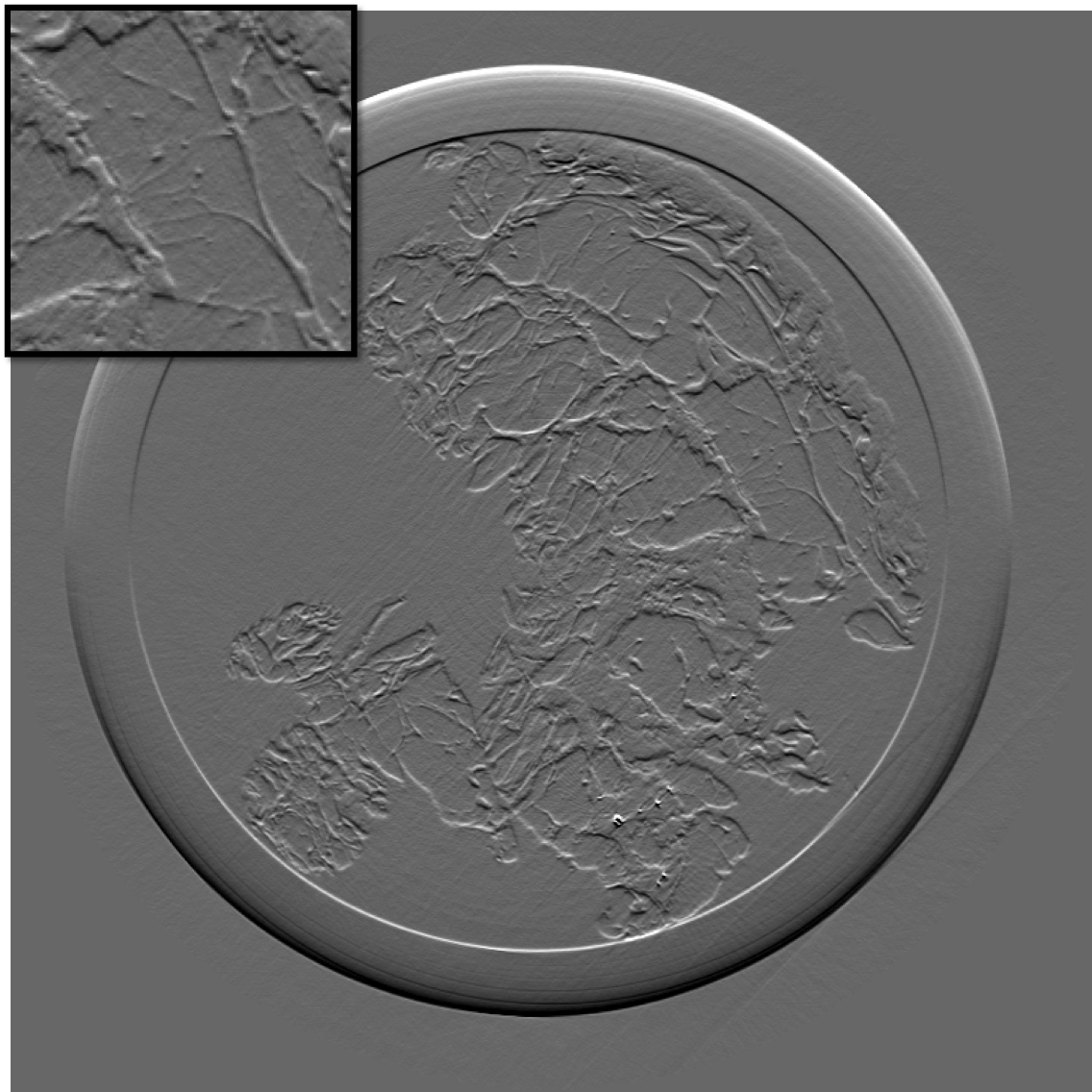}} 
  \hfill{}
   \subfloat{\includegraphics[width=0.245\columnwidth]{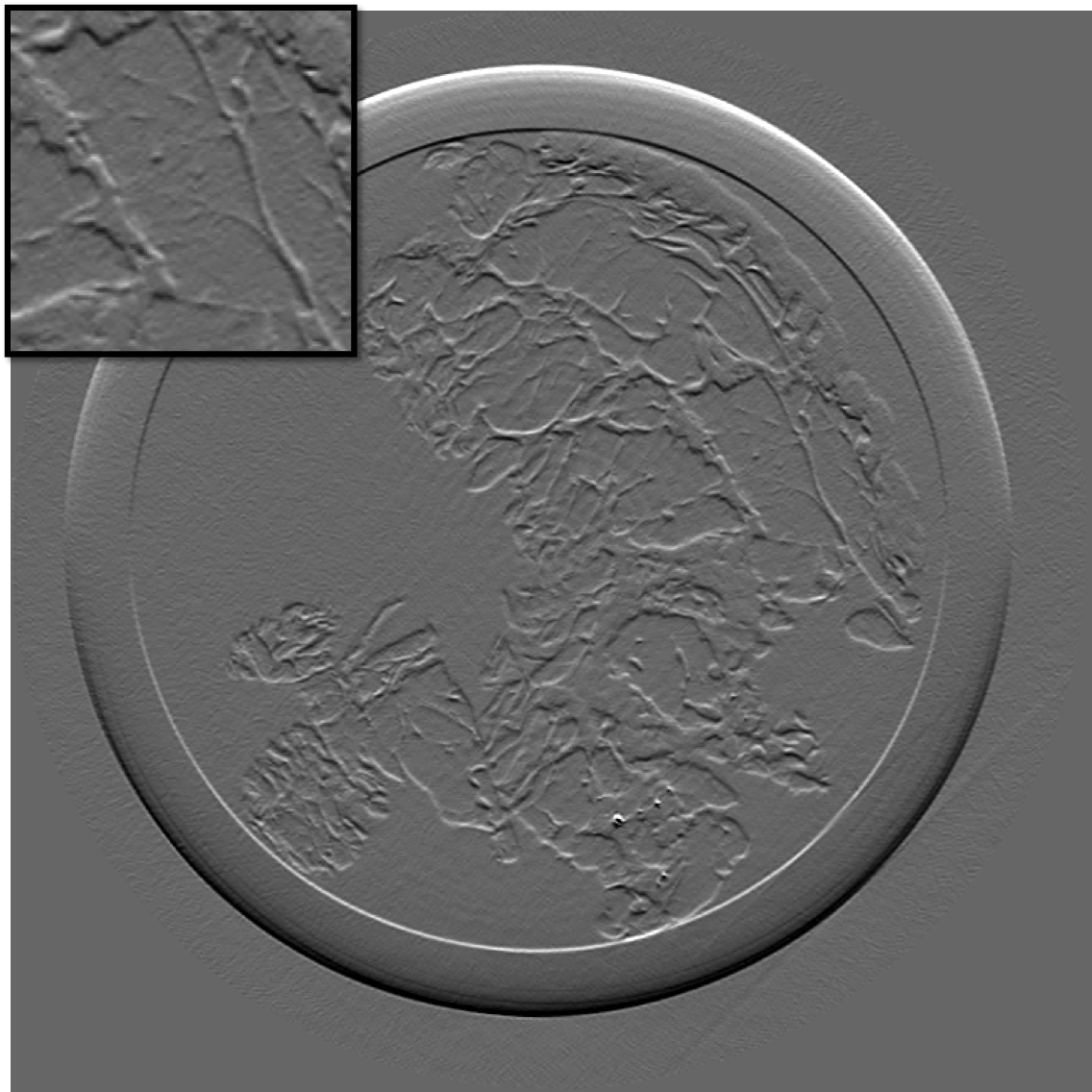}} 
 
 \caption{Phase gradient X(first row) and Y(second row) components reconstruction. Left : FBP with the full set of data. Right: our method with one projection over five, using $\beta=3*10^{-6}$
and $\rho=10$.\label{fig:breastcosfbppp}}
 \end{center}
 
\end{figure}

The sample studied is a 7cm human breast imaged with a pixel size of 100 $\mu$m. The experiment was conducted at the biomedical beam line at the European Synchrotron Radiation Facility (ESRF). The sample was a human breast mastectomy specimen. The study was performed in accordance with the Declaration of Helsinki and was approved by the local ethics committee. A monochromatic X-ray beam with energy of 60 keV was used to image the breast cancer sample. Result of reconstruction is shown in figure \ref{fig:breastcosfbppp}(top). On this image, radiologists could easily identify the skin, fat and glandular tissue.

The training set is obtained from another breast sample imaged by  the same technique but with an high quality reconstruction. We consider a slice image, we apply a Sobel filter to extract the two derivative components and we run the KSVD algorithm.

\begin{figure}
\begin{center}

\includegraphics[width=0.8 \columnwidth]{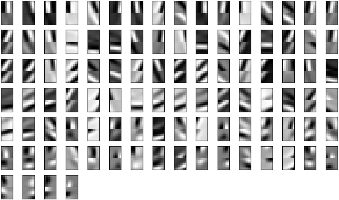}%
\caption{ The vectorial basis of patches learned from a high quality reconstruction of the phase gradient of a human breast. In each patch the upper $7\times 7$ part is the $X$ component while the lower part is the $Y$ component. \label{fig:breast_patches}}

\end{center}
\end{figure}

Figure \ref{fig:breast_patches} shows the patches basis functions that we use to fit both components at the same time. The patches size is $7\times 7$ pixels and each basis function is displayed as a $14x7$ rectangle whose upper $7\times 7$ part is the $X$ component and the lower one the $Y$ component. 

Figure \ref{fig:breastcosfbppp} is the reconstruction for a $765\times 765$ pixels slice, using only 200 projections over the 1000 available. The upper left square is a zoom in the region marked in sub figure \ref{fig:breastcosfbppp}. The right column is the reconstruction with our method for X and Y components, while the left column is reconstructed with standard filtered back-projection using all 1000 available projections. Using our method we can still generate a high quality image with only one fifth of the projections which would be otherwise necessary to generate a high quality reconstruction with the standard FBP method. Visually the difference between the FBP results obtained with full data set and our method with a five-fold reduction of data is barely noticeable. The different borders of structures like skin layers, fatty tissues, and collagen strands are easily identified. The obtained result are very promising and  a systematic evaluation for clinical application is under-way. The radiation dose absorbed by the sample during 200 projections is comparable to that of a standard clinical dual view (2D) mammography (3.5mGy).

%
%
%
%
%
%
%
%

For an eventual future clinical application of the PCI method it is important to investigate which is the acceptable compromise in terms of low dose and sufficient  level of image quality. We need therefore to better explore how the quality of the reconstruction is degraded when we reduce the dose (i.e. number of projections and the acquisition time) further below the standard values. To this end, we performed a reconstruction with only 125 projections and results are shown in the figure III-D for one gradient differential image.

\begin{figure}[htp]
\begin{center}
 \subfloat{\includegraphics[width=0.492\columnwidth]{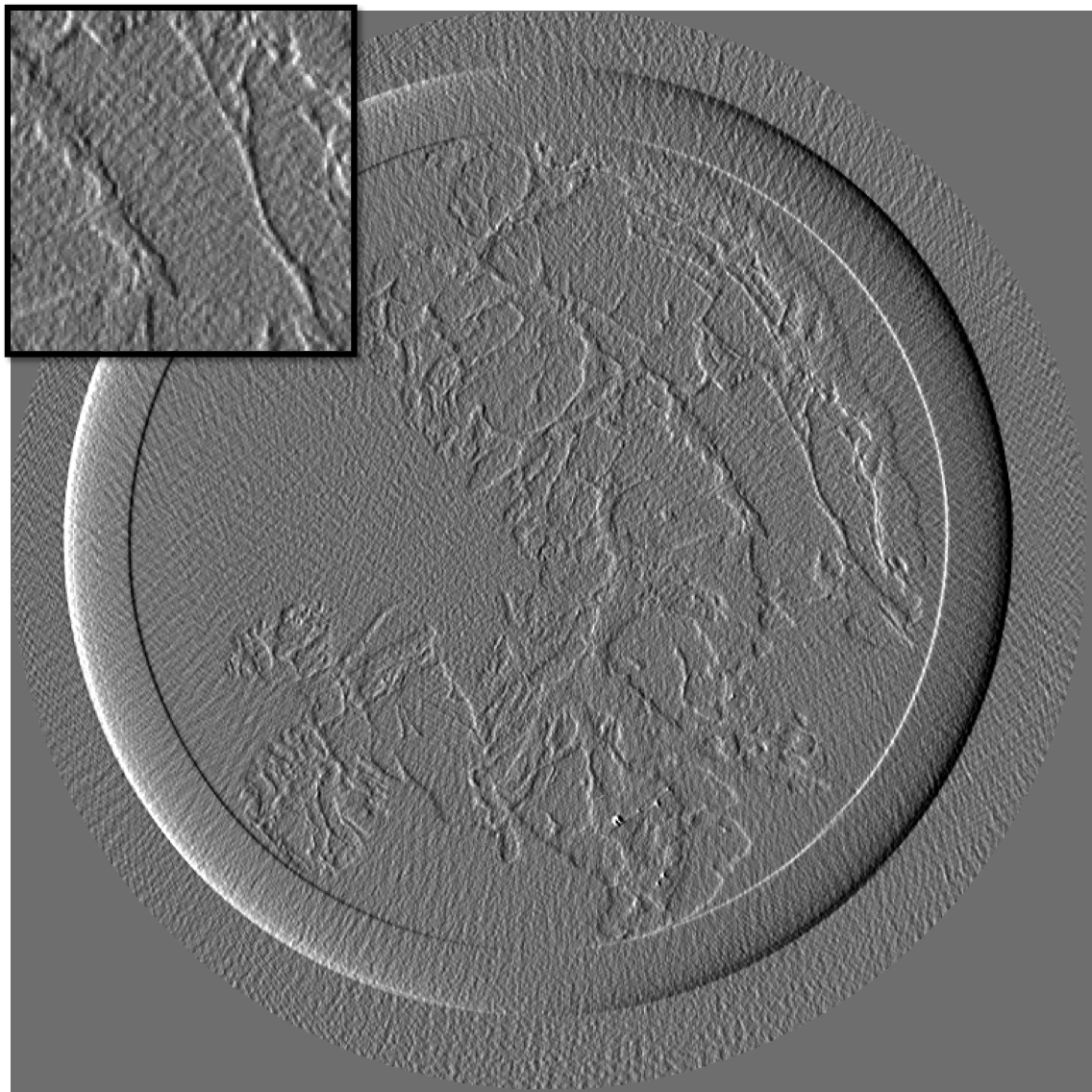}} 
  \hfill{}
   \subfloat{\includegraphics[width=0.492\columnwidth]{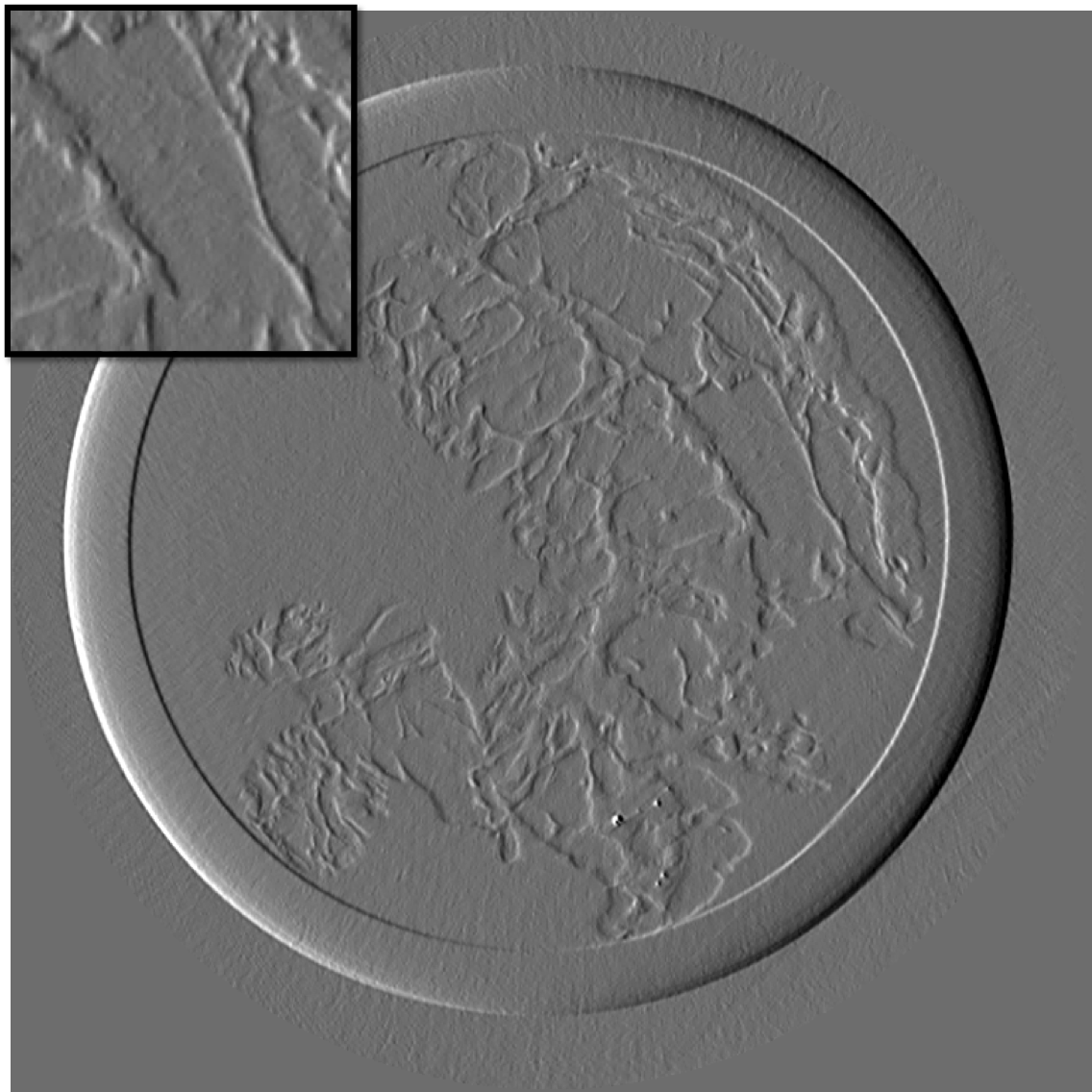}} 
 
 \subfloat{\includegraphics[width=0.492\columnwidth]{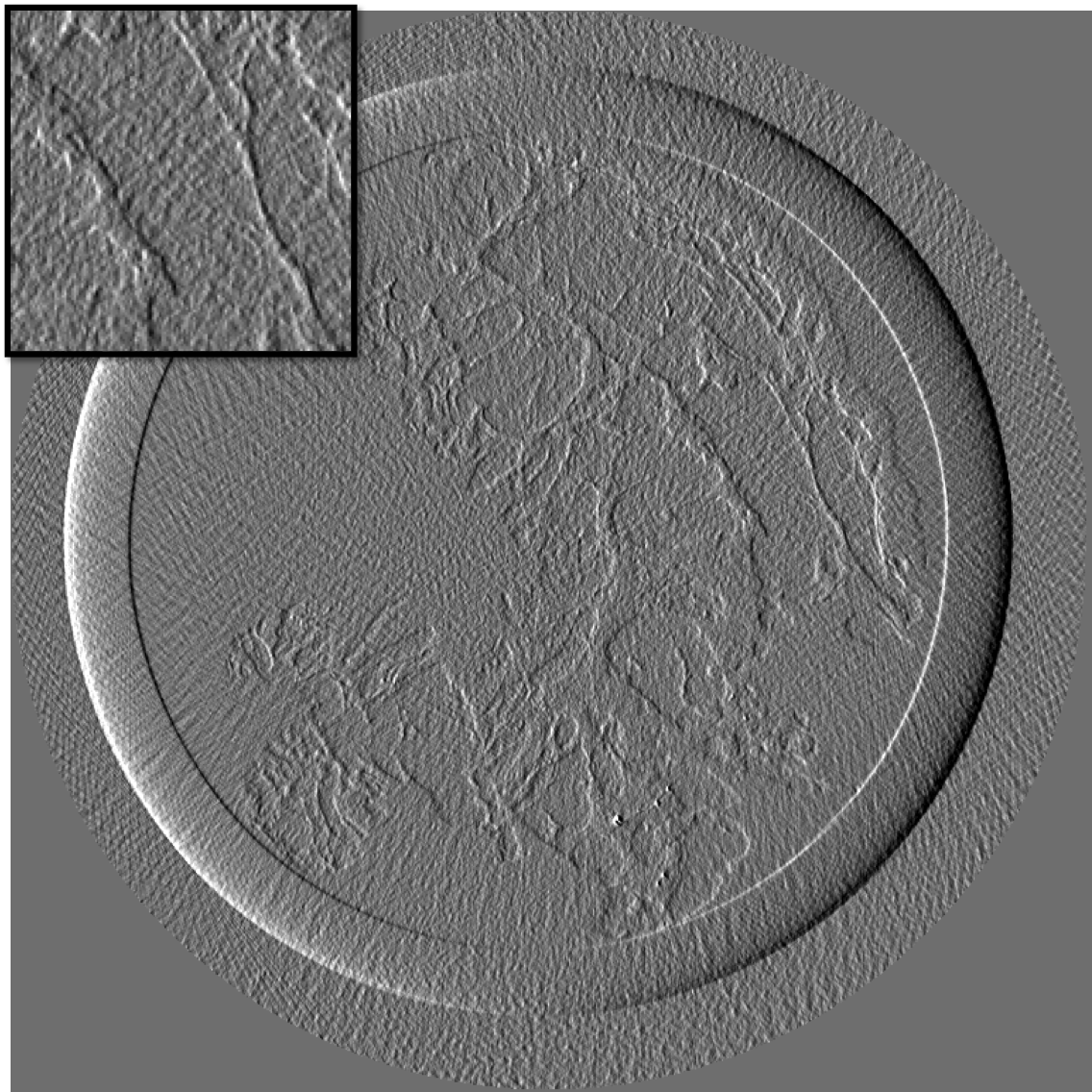}} 
  \hfill{}
   \subfloat{\includegraphics[width=0.492\columnwidth]{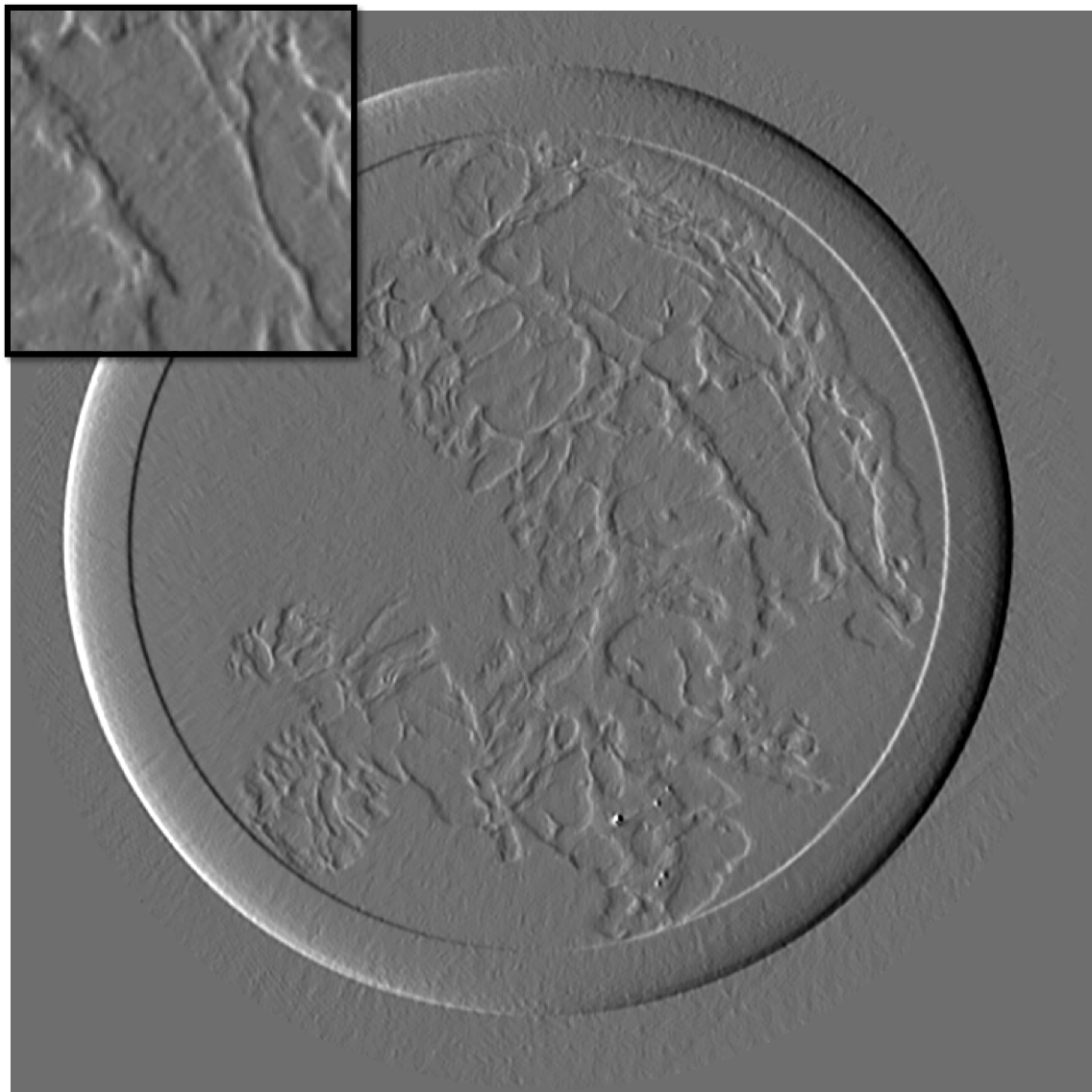}} 
 
 \caption{Phase gradient X reconstruction. Left : FBP with 125 (top) and 100 projections. Right: our method using $\beta=3*10^{-6}$
and $\rho=10$.\label{fig:breastcosfbppp1}}
 \end{center}
 
\end{figure}

For an eventual future clinical application of the PCI method it is important to investigate which is the acceptable compromise in terms of low dose and sufficient  level of image quality. We need therefore to better explore how the quality of the reconstruction is degraded when we reduce the dose (i.e. number of projections and the acquisition time) further below the standard values. To this end, we performed a reconstruction with only 125 projections and results are shown in the figure\ref{fig:breastcosfbppp1}. The first column present the result using our method, the second column is the result of reconstruction using FBP algorithm.

If a slightly higher noise level is tolerable, the method may be used with very few projections and thus applied to the screening and diagnosis of human breast cancers with an even lower radiation dose than conventional dual mammography. The results of our reconstruction show an image quality and a capability of discriminating fine structures that are still clinically acceptable. On the contrary, images produced with the standard FBP reconstruction method are very noisy and not diagnostically satisfactory.
	
%
%
%
%
%
%

\section{CONCLUSION}
We have presented a new convex functional which implements in a mathematically pure form the concept of overlapping-patches-averaging, which was used so-far with a non-convex formalism. The resulting algorithm is robust, efficient, and well adapted to strongly reduce the noise in a natural image. The method was applied also to a medical diagnostic case by considering phase contrast tomographic data of whole cancer bearing human breasts acquired with phase contrast imaging. We demonstrated that it is possible to reduce the deposited dose  in breast CT by a factor 5 compared to the standard algorithm while keeping the same image quality. Although we used this specific example as proof of principle in this study, the method we developed and described can be easily applied to other medical tomography fields. The numerical results have been generated with PyHST \cite{PyHST2,pyhst2_gpl}, tthe ESRF tomography reconstruction code which uses the GPU  implementation of the presented methods.

\section*{Acknowledgment}

We thank Emmanuelle Gouillart and Gael Varoquaux for the interesting discussions and for having pointed us the possibilities of dictionary learning techniques and the examples contained in the scikit-learn package\cite{scikit-learn}. We thank the ESRF for providing the experimental facilities, the ESRF ID17 team for assistance in operating the facilities. We thank Alberto Bravin, Sergei Gasilov, Alberto Mittone for their help during experiments. This work was partially supported by the DFG-Cluster of Excellence Munich-Centre for Advanced Photonics EXE158.

\bibliographystyle{unsrt}
\addcontentsline{toc}{section}{\refname}\bibliography{bib}

\end{document}